\documentclass[12pt]{amsart}
\usepackage{amsfonts}
\usepackage{amsmath,amssymb,amsthm}

\newcommand{\R}{\mathbb R}

\newcommand{\ds}{\displaystyle}
\newtheorem{thm}{Theorem}[section]
\newtheorem{cor}[thm]{Corollary}
\newtheorem{lem}[thm]{Lemma}
\newtheorem{prop}[thm]{Proposition}
\theoremstyle{definition}
\newtheorem{defn}[thm]{Definition}
\theoremstyle{remark}

\begin{document}

\title[ Time-like Weingarten surfaces with real principal curvatures ]
{ Time-like Weingarten surfaces with real principal curvatures in
the three-dimensional Minkowski space and their natural partial
differential equations }%

\author{Vesselka Mihova and Georgi Ganchev}%
\address{Faculty of Mathematics and Informatics, University of Sofia,
J. Bouchier Str. 5, 1164 Sofia, Bulgaria}
\email{mihova@fmi.uni-sofia.bg}%
\address{Bulgarian Academy of Sciences, Institute of Mathematics and Informatics,
Acad. G. Bonchev Str. bl. 8, 1113 Sofia, Bulgaria}%
\email{ganchev@math.bas.bg}

\subjclass[2000]{Primary 53A10, Secondary 53A05}%
\keywords{Time-like W-surfaces in Minkowski space, natural
parameters on time-like W-surfaces in Minkowski space, natural PDE's
of time-like
W-surfaces in Minkowski space.}%

\begin{abstract}

We study time-like surfaces in the three-dimensional Minkowski space with
diagonalizable second fundamental form. On any time-like W-surface we
introduce locally natural principal parameters and prove that
such a surface is determined uniquely (up to motion) by a special invariant
function, which satisfies a natural non-linear partial differential equation.
This result can be interpreted as a solution of the Lund-Regge reduction
problem for time-like W-surfaces with real principal curvatures in Minkowski
space. We apply this theory to linear fractional time-like W-surfaces and
obtain the natural partial differential equations describing them.
\end{abstract}

\maketitle

\section{Introduction}
It has been known to Weingarten \cite{W1, W2}, Eisenhart \cite{E},
Wu \cite{Wu}) that without changing the principal lines on a
Weingarten surface in  Euclidean space, one can find geometric coordinates in which the
coefficients of the metric are expressed by the principal curvatures
(or principal radii of curvature).

The geometric parameters on Weingarten surfaces were used in
\cite{Wu} to find the classes of Weingarten surfaces yielding
``geometric" $\mathfrak{so}(3)$-scattering systems (real or complex)
for the partial differential equations, describing these surfaces.

We have shown that the Weingarten surfaces in Euclidean space \cite{GM1, GM2}
and space-like surfaces in Minkowski space \cite{GM3} admit geometrically
determined principal parameters (\emph{natural principal parameters}),
which have the following property: all invariant functions on W-surfaces
can be expressed in terms of one function $\nu$, which satisfies one \emph{natural}
partial differential equation. The Bonnet type fundamental theorem states that
any solution to the natural partial differential equation determines a W-surface
uniquely up to motion. Thus the description of any class of W-surfaces (determined
by a given Weingarten relation) is equivalent to the study of the solution space
of their natural PDE. This solves the Lund-Regge reduction problem \cite{LR}
for W-surfaces in Euclidean space and space-like W-surfaces in Minkowski space.

The relationship between the solutions of certain types of partial differential equations
and the determination of various kinds of surfaces of constant curvature has generated
many results which have applications to the areas of both pure and applied mathematics.
This includes the determination of surfaces of either constant mean curvature or
Gaussian curvature. It has long been known that there is a connection between surfaces
of negative constant Gaussian curvature in Euclidean $\mathbb{R}^3$ and the sine-Gordon equation.
The fundamental equations of surface theory are found to yield a type of geometrically
based Lax pair. For instance, given a particular solution of the sinh-Laplace equation,
this Lax pair can be integrated to determine the three fundamental vector fields
related to the surface.
These are also used to determine the coordinate vector field of the surface.

Further results are obtained based on the fundamental equations of surface theory,
and it is shown how specific solutions of this sinh-Laplace equation can be used
to obtain the coordinates of a surface in either Minkowski $\mathbb{R}^3_1$
or Euclidean $\mathbb{R}^3$ space \cite{Hu1, Hu2}.

In \cite{Br} Bracken introduces some fundamental concepts and equations pertaining
to the theory
of surfaces in three-space, and, in particular, studies a class of sinh-Laplace equation
which has the form
$\Delta u=\pm \sinh u.$

In this paper we study time-like surfaces with real principal curvatures
in the three dimensional Minkowski space $\R^3_1$.

A time-like surface $\mathcal M$ with real principal curvatures $\nu_1$
and $\nu_2$ is a Weingarten surface (W-surface) \cite{W1, W2} if there exists a
function $\nu$ on $\mathcal M$ and two functions (Weingarten functions) $f, \, g$ of one
variable, such that
$$\nu_1=f(\nu), \quad \nu_2=g(\nu).$$

A basic property of W-surfaces in Euclidean space is the following
theorem of Lie \cite{Lie}:

\emph{The lines of curvature of any W-surface can be found in quadratures.}

This remarkable property is also valid for space-like and time-like W-surfaces in
Minkowski space.

We use four invariant functions (two principal normal curvatures $\nu_1,\,\nu_2$
and two principal geodesic curvatures $\gamma_1,\, \gamma_2$) and divide time-like
W-surfaces into two classes with respect to these invariants:

(1) the class of \emph{strongly regular} time-like surfaces defined by
$$(\nu_1-\nu_2)\,\gamma_1 \,\gamma_2\neq 0;$$

(2) the class of time-like surfaces defined by
$$\gamma_1= 0, \quad (\nu_1-\nu_2)\, \gamma_2\neq 0.$$

The basic tool to investigate the relation between time-like surfaces and
the partial differential equations describing them, is Theorem 2.1. This theorem
is a reformulation of the fundamental Bonnet theorem for the class of strongly
regular time-like surfaces in terms of the four invariant functions.
Further, we apply this theorem to time-like W-surfaces.

In Section 3 we prove (Proposition 3.3) that any time-like W-surface admits
locally special principal parameters (\emph{natural principal parameters}).

Theorem 3.6 is the basic theorem for time-like W-surfaces of type (1):

\emph{Any strongly regular time-like W-surface is determined uniquely up to motion by
the functions $f$, $g$ and the function $\nu$, satisfying the natural PDE $(3.3)$.}

Theorem 3.7 is the baic theorem for time-like Weingarten surfaces of
type (2):

\emph{Any time-like W-surface with $\gamma_1=0$ is determined
uniquely up to motion by the functions $f, g$ and the function
$\nu$, satisfying the natural ODE $(3.8)$.}

In natural principal parameters the four basic invariant functions, which determine
time-like W-surfaces uniquely up to motions in $\R^3_1$, are expressed by a single
function, and the system of Gauss-Codazzi equations reduces to a single partial
differential equation (the Gauss equation). Thus, the number of the four invariant
functions, which determine time-like W-surfaces, reduces to one
invariant function, and the number of Gauss-Codazzi equations reduces to one
\emph{natural} PDE. This result gives a solution to the Lund-Regge reduction
problem \cite{LR} for the time-like W-surfaces in $\R^3_1$. The Lund-Regge
reduction problem has been analyzed and discussed from several view points in
the paper of Sym \cite{Sym}.

In Proposition 4.1 we prove that

{\it The natural principal parameters of a given time-like W-surface $\mathcal{M}$
are natural principal parameters for all parallel time-like surfaces
$\overline{{\mathcal{M}}}(a),\; a={\rm const} \neq 0$ of $\mathcal{M}$.}

Theorem 4.2 states that (cf \cite{GM2,GM3}):

{\it The natural PDE of a given time-like W-surface $\mathcal M$ is
the natural PDE of any parallel
time-like surface $ \overline{{\mathcal M}}(a),\; a={\rm const}\neq 0$, of $\mathcal M$.}

In \cite{Mil1, Mil3} Milnor studies surface theory in Euclidean and Minkowski space,
considering harmonic maps and various relations between the Gauss curvature $K,$
the mean curvature $H$ and the curvature
$H'=\ds{\frac{\nu_1-\nu_2}{2}}$.
In \cite {Mil2, GM2} is proved  that any  surface in $\R^3_1$,
whose Gauss curvature $K$ and mean curvature $H$  satisfy the linear relation
$$\delta K = \alpha H + \gamma, \quad \alpha, \gamma, \delta -
{\rm constants}; \quad \alpha^2 + 4 \gamma\delta \neq 0,\leqno(1.1)$$
is parallel to a surface, satisfying one of the following conditions:
$H=0$, $K=1$ or $K=-1$.

There arises the following question: what are the natural PDE's describing the surfaces,
whose curvatures satisfy the relation (1.1)?

Since any time-like surface $\mathcal M$, whose invariants $K$ and $H$ satisfy the linear
relation (1.1), is (locally) parallel to one of the following three types of basic
surfaces: a surface with $H=0$; a surface with $K=1$; a surface
with $K=-1$, from Theorem 4.2 it follows that

{\it Up to similarity, the time-like surfaces, whose curvatures satisfy
the linear relation $(1.1)$,
are described by the natural PDE's of the basic surfaces.}

A. Ribaucour \cite{Rib} has proved that
{\it a necessary condition for the curvature lines
of the first and second focal surfaces of $\mathcal{M}$ to  correspond to each other resp. to
a conjugate parametric lines on $\mathcal{M}$ is $\rho_1-\rho_2={\rm const}$ resp.
$\rho_1\,\rho_2={\rm const}$}.

Von Lilienthal \cite{vL1} (cf \cite{vL2, Bi1, Bi2, E}) has proved in $\mathbb{R}^3$
that a surface with a relation
$\displaystyle{\rho_1-\rho_2=\frac{1}{R},\; R={\rm const}\neq 0},$
between its principal radii of curvature $\displaystyle{\rho_1=\frac{1}{\nu_1}}$ and
$\displaystyle{\rho_2=\frac{1}{\nu_2}}$
has  first and second focal surfaces $\widetilde {\mathcal M}$ of constant Gauss
curvature $-R^2$ and vice versa.
The involute surfaces  $ \overline{{\mathcal M}}(a), \, a \in\mathbb{R}$ of $ \widetilde{{\mathcal M}}$
are parallel surfaces of $ \mathcal M$ with the property $\rho_1-\rho_2={\rm const}$.
This implies that the family $\overline{{\mathcal M}}(a)$ are integrable surfaces
as a consequence of the integrability of $\widetilde{\mathcal M}$.
The curvatures of the above surfaces $ \mathcal M$ satisfy the relation
$K=\beta\,H',\; \beta={\rm const}\neq 0$.

In $\mathbb{R}^3_1$ one can prove in a similar way the corresponding property:
The first focal surface of a time-like surface with $K=\beta\,H',\, \beta\neq 0,$ is
space-like  of constant Gauss curvature $\beta^2/4$, and its second focal surface is
time-like of constant Gauss curvature $-\beta^2/4$.

Obviously the time-like surfaces with $K=\beta\,H',\; \beta={\rm const}\neq 0,$ are not
included in the class characterized by (1.1).

These surfaces belong to the classes of time-like W-surfaces, defined by the
following more general linear relation
$$\delta K = \alpha H + \beta H' + \gamma, \quad \alpha, \beta, \gamma, \delta -
{\rm constants}; \quad \alpha^2-\beta^2 + 4 \gamma\delta \neq 0\leqno(1.2)$$
 between the Gauss curvature $K$, the mean curvature $H$ and the curvature $H'$.
We denote this class by $\mathfrak{K}$.

We show that the class $\mathfrak{K}$ is  the class
of linear fractional time-like W-surfaces with respect to the principal curvatures
(cf \cite{GM2,GM3}).
Furthermore, if $\mathcal M$ is a time-like surface in $\mathfrak{K}$,
then its parallel surfaces $ \overline{{\mathcal M}}(a),\, a={\rm const},$
belong to $\mathfrak{K}$ too.

In the main Theorem 5.3 in this paper we determine ten basic relations
with respect to the constants in (1.2) and each
of them generates a {\it basic subclass of surfaces} of $\mathfrak{K}$.
Any time-like surface $\mathcal M$, whose invariants $K$, $H$ and $H'$ satisfy the linear
relation (1.2) is (locally) parallel to one of these basic surfaces.

In \cite{Hu2} Hu has cleared up the relationship between the PDE's
$$\begin{array}{ll}
\alpha_{uu}-\alpha_{vv}=\pm\sin \alpha &(\sin{\rm -Gordon \; PDE}), \\
[2mm]
\alpha_{uu}-\alpha_{vv}=\pm\sinh \alpha &(\sinh{\rm -Gordon \; PDE}), \\
[2mm]
\alpha_{uu}+\alpha_{vv}=\pm\sin \alpha &(\sin{\rm -Laplace \; PDE}), \\
[2mm]
\alpha_{uu}+\alpha_{vv}=\pm\sinh \alpha &(\sinh{\rm -Laplace \; PDE})
\end{array}$$
and the construction of various kinds of surfaces of constant curvature in $\mathbb{R}^3$
or $\mathbb{R}^3_1$.

In \cite{Hu3} by using Darboux transformations, from a known solution to
the $\sinh$-Laplace (resp. $\sin$-Laplace) equation  have been obtained explicitly
new solutions to the $\sin$-Laplace (resp. $\sinh$-Laplace) equation.

Time-like surfaces with positive Gauss curvature and imaginary principal curvatures
have been constructed in \cite{Gu}.

\vskip 2mm
It is essential to note that the natural PDE's of the time-like W-surfaces
from the class $\mathfrak{K}$ are expressed in the form $ \delta
\lambda = f(\lambda)\,$, \;where $\delta$ is one of the operators
(cf \cite{GM2, GM3}):
$$\Delta \lambda:=\lambda_{xx}+\lambda_{yy},  \qquad
\bar \Delta \lambda:=\lambda_{xx}-\lambda_{yy};$$
$$\Delta^* \lambda:=\lambda_{xx}+(\lambda^{-1})_{yy},  \qquad
\bar \Delta^* \lambda:=\lambda_{xx}-(\lambda^{-1})_{yy}.$$
\vskip 3mm

\section{Preliminaries}

Let $\R^3_1$ be the three dimensional Minkowski space with the standard flat
metric $\langle \, , \, \rangle$ of signature $(2,1)$. We assume that the
following orthonormal coordinate system $Oe_1e_2e_3: \; e_1^2=e_2^2=-e_3^2=1, \;
\langle e_i, \, e_j \rangle = 0, \, i\neq j$ is fixed and gives the orientation
of the space.

Let $\mathcal M: \, z=z(u,v), \; (u,v) \in {\mathcal D}$ be a
time-like surface in the three dimensional Minkowski space ${\R}_1^3$ and
$\nabla$ be the flat Levi-Civita connection of the metric $\langle \, , \, \rangle$.
The unit normal vector field to ${\mathcal M}$ is denoted by $l$ and
$E, F, G; \; L, M, N$ stand for the coefficients of the first and the second
fundamental forms, respectively. Then we have
$$E=z_u^2<0, \quad F=z_u\,z_v,\quad G=z_v^2>0, \quad EG-F^2<0, \quad l^2=1.$$

The coefficients of the second fundamental form are given as follows:
$$L= l\,z_{uu}=-l_u\,z_u, \quad M=l\,z_{uv}=-l_u\,z_v
=-l_v\,z_u, \quad N= l\,z_{vv}=-l_v\,z_v.$$

The linear Weingarten map $\gamma$ is determined by the conditions
$$\gamma(z_u)=l_u, \quad \gamma(z_v)=l_v.$$
Then the mean curvature $H$ and the Gauss curvature $K$ of $\mathcal M$ are given in the standard way
$$H=-\frac{1}{2}\, {\rm tr \, \gamma}, \quad K={\rm det}\, \gamma.$$

While the Weingarten map of a space-like surface satisfies the
inequality $H^2-K \geq 0$ and is always diagonalizable, the Weingarten map on a
time-like surface can satisfy the inequalities $H^2-K \geq 0$ or $H^2-K < 0$.

Throughout the whole paper we deal with time-like surfaces satisfying the
inequality
$\;H^2-K \geq 0,$ i.e. time-like surfaces with real principal curvatures.

We suppose that the surfaces under consideration are free of points with $H^2-K=0$,
i.e. satisfy the strong inequality
$$H^2-K > 0\leqno(2.1)$$
and denote by $H'$ the invariant curvature
$$H'=\sqrt{H^2-K}.$$

Under the above condition the theory of time-like surfaces can be developed in a way
similar to the theory of surfaces in Euclidean space or space-like surfaces in
Minkowski space.

Time-like surfaces satisfying the condition (2.1) can be locally parameterized by
principal parameters. Further we assume that the parametric net is principal, i.e.
$$F(u,v)=M(u,v)=0, \quad (u,v) \in \mathcal D.$$

Then the principal curvatures $\nu_1, \nu_2$ and the principal geodesic curvatures
(geodesic curvatures of the principal lines) $\gamma_1, \gamma_2$ are given by
$$\nu_1=\frac{L}{E}, \quad \nu_2=\frac{N}{G}; \qquad
\gamma_1=\frac{E_v}{2E\sqrt G}, \quad \gamma_2= \frac{-G_u}{2G\sqrt {-E}}, \leqno(2.2)$$
and $\nu_1$, $\nu_2$ satisfy the Rodrigues' formulas:
$$l_u=-\nu_1\,z_u, \quad l_v=-\nu_2\,z_v.$$

We consider the tangential frame field $\{X, Y\}$ determined by
$$X:=\frac{z_u}{\sqrt {-E}}\,, \qquad Y:=\frac{z_v}{\sqrt G}$$
and suppose that the moving frame field $XYl$ is positive oriented.

The following Frenet type formulas for the frame field $XYl$ are valid

$$\begin{tabular}{ll}
$\left|\begin{array}{llccc}
\nabla_{X} \,X & = &  &\gamma_1 \,Y - \nu_1 \, l,  &\\
[2mm]
\nabla_{X} Y & = \; \; \, \gamma_1 \, X, & & \\
[2mm]
\nabla_{X} \, l & = - \nu_1 \, X, & & &
\end{array}\right.$ &
\quad $\left|\begin{array}{llccc}
\nabla_{Y} \,X & = & & -\gamma_2 \, Y, &\\
[2mm]
\nabla_{Y} Y & = \; \; \, -\gamma_2 \, X & &  & + \nu_2 \, l,\\
[2mm]
\nabla_{Y}\, l & = &  & - \nu_2 \, Y. &
\end{array}\right.$
\end{tabular}\leqno(2.3)$$
\vskip 2mm

The Codazzi equations have the form
$$\gamma_1=\frac{-Y(\nu_1)}{\nu_1-\nu_2}= \frac{-(\nu_1)_v}{\sqrt G(\nu_1-\nu_2)},
\qquad
\gamma_2=\frac{-X(\nu_2)}{\nu_1-\nu_2}=
\frac{-(\nu_2)_u}{\sqrt {-E}\,(\nu_1-\nu_2)}, \leqno(2.4)$$
and the Gauss equation can be written as follows:
$$X(\gamma_2)+Y(\gamma_1)+ \gamma_1^2-\gamma_2^2 =-\nu_1\nu_2 = -K,$$
or
$$\frac{(\gamma_2)_u}{\sqrt{-E}}+\frac{(\gamma_1)_v}{\sqrt G}+
\gamma_1^2-\gamma_2^2=-\nu_1\nu_2=-K.\leqno(2.5)$$

A time-like surface ${\mathcal M}: \; z=z(u,v), \; (u,v)\in \mathcal D$
parameterized by principal parameters is said to be \emph{strongly regular}
if (cf \cite{GM1, GM2, GM3})
$$(\nu_1(u,v)-\nu_2(u,v))\gamma_1(u,v)\gamma_2(u,v) \neq 0, \quad (u,v) \in \mathcal D.$$

The Codazzi equations (2.4) imply that
$$\gamma_1 \gamma_2 \neq 0 \; \iff \; (\nu_1)_v (\nu_2)_u \neq 0.$$
Because of (2.4) the formulas
$$\sqrt {-E}=\frac{-(\nu_2)_u}{\gamma_2(\nu_1-\nu_2)} >0, \quad
\sqrt G=\frac{-(\nu_1)_v}{\gamma_1 (\nu_1-\nu_2)}>0 \, \leqno(2.6)$$
are valid on strongly regular time-like surfaces.

Taking into account (2.6) for strongly regular time-like surfaces formulas (2.3) become

$$\left|\begin{array}{l}
X_u  = \displaystyle{-\frac{\gamma_1 \, (\nu_2)_u}
{\gamma_2(\nu_1-\nu_2)}\, Y}
+ \displaystyle{\frac{\nu_1 \, (\nu_2)_u}{\gamma_2
(\nu_1-\nu_2)}\, l},\;
Y_u  =  - \displaystyle{\frac{\gamma_1 \, (\nu_2)_u}
{\gamma_2(\nu_1-\nu_2)}\, X},\;
l_u  =   \displaystyle{\frac{\nu_1 \, (\nu_2)_u}
{\gamma_2(\nu_1-\nu_2)}}\, X;  \\
[6mm]
X_v =\displaystyle{\frac{\gamma_2 \, (\nu_1)_v}
{\gamma_1(\nu_1-\nu_2)}\, Y,}\;
Y_v  = \displaystyle{\frac{\gamma_2 \, (\nu_1)_v}
{\gamma_1(\nu_1-\nu_2)}}\, X -\displaystyle{\frac{\nu_2 \,(\nu_1)_v}{\gamma_1
(\nu_1-\nu_2)}}\, l,\;
l_v  = \displaystyle{\frac{\nu_2 \, (\nu_1)_v}
{\gamma_1(\nu_1-\nu_2)}}\,Y.
\end{array}\right.\leqno (2.7)$$

Finding the compatibility conditions for the systen (2.7), we reformulate
the fundamental Bonnet theorem for strongly regular time-like surfaces
in terms of the invariants of the surface.

\begin{thm}\label{T:2.1}
Given four functions $\nu_1(u,v), \, \nu_2(u,v), \, \gamma_1(u,v),
\, \gamma_2(u,v)$ defined in a neighborhood $\mathcal D$ of $(u_0, v_0)$
and satisfying the conditions
$$\begin{array}{ll}
1) & (\nu_1-\nu_2)\,\gamma_1 \, (\nu_1)_v <0, \quad
(\nu_1-\nu_2)\,\gamma_2 \, (\nu_2)_u < 0, \\
[4mm]
2.1) & \displaystyle{\left(\ln\frac{(\nu_1)_v}{\gamma_1}\right)_u=
\frac{(\nu_1)_u}{\nu_1-\nu_2},}
\qquad
\displaystyle{\left(\ln\frac{(\nu_2)_u}{\gamma_2}\right)_v=
-\frac{(\nu_2)_v}{\nu_1-\nu_2},}\\
[5mm]
2.2) & \displaystyle{\frac{\nu_1-\nu_2}{2}\left(\frac{(\gamma_2^2)_u}{(\nu_2)_u}+
\frac{(\gamma_1^2)_v}{(\nu_1)_v}\right)- (\gamma_1^2-\gamma_2^2)=\nu_1\nu_2}.
\end{array}$$
Let $z_0X_0Y_0l_0$ be an initial positive oriented orthonormal frame.

Then there exists a unique strongly regular time-like surface
${\mathcal M}: \; z=z(u,v), \; (u,v) \in \mathcal D_0 \;
((u_0, v_0) \in \mathcal D_0 \subset \mathcal D)$ with prescribed invariants
$\nu_1, \, \nu_2, \, \gamma_1, \, \gamma_2$ such that
$$z(u_0, v_0)=z_0, \; X(u_0, v_0)=X_0, \; Y(u_0, v_0)=Y_0, \; l(u_0, v_0)=l_0.$$
\end{thm}

Formulas (2.3) imply explicit expressions for the curvature and
the torsion of any principal line on the time-like surface ${\mathcal M}$.

Let $c_1: z=z(s), \;\mathcal M \in J$ be a line from the family $\mathcal F_1$
($v={\rm const}$) parameterized by a natural parameter and
$\kappa_1, \, \tau_1$ be its curvature and torsion,
respectively.

Since $c_1$ is an integral line of the unit time-like vector field $X$, then
$$z'=X, \quad z''=\nabla_XX=\gamma_1\,Y-\nu_1\,l,$$
$$z'''=\nabla_X\nabla_XX=-X(\nu_1)\,l+X(\gamma_1)\,Y +(\nu_1^2+\gamma_1^2)\,X,$$
$$\kappa_1^2= \nu_1^2+\gamma_1^2.$$
We use the formula
$$\tau=\frac{z'z''z'''}{z''^2}.$$

Since $\nu_1^2+\gamma_1^2>0$ along $c_1$, we find
$$\tau_1=\frac{\nu_1\, X(\gamma_1)-\gamma_1\,X(\nu_1)}
{\nu_1^2+\gamma_1^2}= \frac{
\nu_1^2}{\kappa_1^2}\,X\left(\frac{\gamma_1}{\nu_1} \right).$$
Denoting $\sin\theta_1= \frac{\gamma_1}{\kappa_1}$ and
$\cos\theta_1= \frac{\nu_1}{\kappa_1}$, we obtain
$$\tau_1=X(\theta_1).$$

For the lines $c_2$ of the family $\mathcal F_2$ we obtain in a
similar way the corresponding formulas
$$z'=Y, \quad z''=\nabla_YY=-\gamma_2\,X+\nu_2\,l,$$
$$z'''=\nabla_Y\nabla_YY=Y(\nu_2)\,l-Y(\gamma_2)\,X +(\gamma_2^2-\nu_2^2)\,Y,$$
$$\kappa_2^2=\varepsilon_2\,z''^2=\varepsilon_2\,(\nu_2^2-\gamma_2^2),
\quad \varepsilon_2={\rm sign}\,z''^2,$$
and in the case $z''^2\neq 0$
$$\tau_2=\varepsilon_2\,\frac{\gamma_2\, Y(\nu_2)-\nu_2\,Y(\gamma_2)}
{\kappa_2^2}=-\varepsilon_2\,\frac{\nu_2^2}{\kappa_2^2}\,Y\left(\frac{\gamma_2}{\nu_2}\right)
.$$

\vskip 4mm
\section{Natural principal parameters on time-like Weingarten surfaces}

In this section we consider diagonalizable  time-like Weingarten surfaces.
For the sake of symmetry with respect to the principal curvatures
$\nu_1$ and $\nu_2$ we use the following characterization of
 time-like Weingarten surfaces:

\emph{A diagonalizable time-like surface $\mathcal M:\; z=z(u,v),\; (u,v)\in
\mathcal{D}$ is Weingarten if there exist two real differentiable
functions $f(\nu), \; g(\nu), \; f(\nu)-g(\nu)\neq 0,\;
f'(\nu)g'(\nu)\neq 0, \; \nu \in \mathcal{I}\subseteq{\R}$ such that
the principal curvatures of $\mathcal M$ at every point are given by
$\nu_1=f(\nu), \; \nu_2=g(\nu), \; \nu=\nu(u,v), \;
(u,v) \in \mathcal D$.} \vskip 2mm

The next statement gives a property of  time-like Weingarten
surfaces, which allows us to introduce special principal parameters
on such surfaces.

\begin{lem} Let $\mathcal M:\; z=z(u,v),\; (u,v)\in \mathcal{D}$ be a
diagonalizable time-like Weingarten surface parameterized with principal parameters.
Then the function
$$\lambda = \sqrt {-E} \exp\left(\int \frac{f'd\nu}{f-g}\right)$$
does not depend on $v$, while the function
$$\mu = \sqrt G \exp \left(\int \frac{g'd\nu}{g-f}\right)$$
does not depend on $u$.
\end{lem}

{\it Proof} : Taking into account (2.4) and (2.2), we find
$$\gamma_1=\frac{-f'(\nu) Y(\nu)}{f(\nu)-g(\nu)}=Y(\ln \sqrt {-E}), \qquad
\gamma_2=\frac{-g'(\nu) \,X(\nu)}{f(\nu)-g(\nu)}= -X(\ln \sqrt G),$$
which imply that
$$Y\left(\int\frac{f'(\nu)\,d \nu}{f(\nu)-g(\nu))} + \ln \sqrt {-E}\right) = 0,
\qquad X\left(\int\frac{g'(\nu)\, d \nu}{g(\nu)-f(\nu)} + \ln \sqrt G\right) = 0.$$
The last equalities mean that $\lambda_v=0$ and
$\mu_u=0$. \qed
\vskip 2mm

We define special principal parameters on a  time-like
Weingarten surface as follows:

\begin{defn}
Let $\mathcal M:\; z=z(u,v),\; (u,v)\in \mathcal{D}$ be a diagonalizable time-like
Weingarten surface parameterized with principal parameters. The
parameters $(u, v)$ are said to be \emph{natural} principal, if
the functions $\lambda(u)$ and $\mu(v)$ from Lemma 3.1 are
constants.
\end{defn}

\begin{prop}\label {P:5.4}
Any diagonalizable time-like Weingarten surface admits locally natural principal parameters.
\end{prop}

\emph{Proof}: Let $\mathcal M:\; z=z(u,v),\; (u,v)\in \mathcal{D}$ be a
 time-like Weingarten surface in the Minkowski space $\mathbb{R}^3_1$,
parameterized with principal parameters. Then $\nu_1=f(\nu), \; \nu_2=g(\nu), \; \nu=\nu(u,v)$
for some differentiable functions $f$, $g$ and $\nu$ satisfying the
conditions
$\;(f(\nu)-g(\nu))\, f'(\nu)\,g'(\nu)\neq 0, \; (u,v) \in \mathcal D.$

Let $\mathfrak{a}={\rm const} \neq 0, \; \mathfrak{b}={\rm const} \neq 0$, \; $(u_0, v_0)
\in \mathcal D$ and $\nu_0=\nu(u_0,v_0)$. We change the parameters
$(u, v)\in \mathcal{D}$ with $(\bar u, \bar v)\in \bar{\mathcal{D}}$
by the formulas
$$\begin{array}{l}
\displaystyle{\bar u=\mathfrak{a}\int_{u_0}^u \sqrt {-E}
\exp\left(\int_{\nu_0}^{\nu} \frac{f'd\nu}{f-g}\right)}
\, du\, +\overline{u}_0 , \quad \bar u_0 = {\rm const},\\
[4mm] \displaystyle{\bar v=\mathfrak{b}\int_{v_0}^v\sqrt G \exp
\left(\int_{\nu_0}^{\nu} \frac{g'd\nu}{g-f}\right)\,
dv}\,+\overline{v}_0, \quad \bar v_0={\rm const}.
\end{array}$$
According to Lemma 3.1 it follows that $(\bar u, \bar v)$
are again principal parameters and
$$\bar E=-\frac{1}{\mathfrak{a}^2}\,\exp \left(-2\int_{\nu_0}^{\nu} \frac{f'd\nu}{f-g}\right),
\quad \bar G=\frac{1}{\mathfrak{b}^2}\,\exp \left(-2\int_{\nu_0}^{\nu}
\frac{g'd\nu}{g-f}\right).\leqno(3.1)$$

Then for the functions from Lemma 3.1  we find
$$ \lambda(\bar u) =|\mathfrak{a}|^{-1}, \quad  \mu(\bar v) =|\mathfrak{b}|^{-1}.$$
Furthermore $\;\mathfrak{a}^2\,\bar E(u_0,v_0)=-1, \; \mathfrak{b}^2\, \bar G(u_0,v_0)=1$.
{\qed} \vskip 2mm

We assume now that the considered  time-like Weingarten
surface $\mathcal M:\, z=z(u,v)$, $(u,v) \in \mathcal D$ is parameterized with
natural principal parameters $(u, v)$. It follows from the above
proposition that the coefficients $E$ and $G$  (consequently $L$ and
$N$) are expressed by the invariants of the surface.

As an immediate consequence from Proposition \ref{P:5.4} we get

\begin{cor}
Let $\mathcal M$ be a  time-like Weingarten surface parameterized by
natural principal parameters $(u, v)$. Then any natural
principal parameters $(\tilde u, \tilde v)$ on $\mathcal M$ are determined by
$(u, v)$ up to an affine transformation of the type
$$ \tilde u = a_{11} \, u + b_1,\quad
\tilde v = a_{22} \, v + b_2, \quad
a_{11}a_{22} \neq 0,$$
or of the type
$$\tilde u = a_{12} \, v  + c_1,\quad
\tilde v = a_{21} \, u +c_2, \quad
a_{12}a_{21}\neq 0,$$
where $a_{ij}, \, b_i, \, c_i; \, i,j=1,2 $ are constants.
\end{cor}

Next we give a simple criterion principal parameters to be
natural.

\begin{prop}\label{P:7.2}
Let a  time-like Weingarten surface $\mathcal M:\; z=z(u,v),\; (u,v)\in
\mathcal{D}$ be parameterized with principal parameters. Then
$(u,v)$ are natural principal if and only if
$$\sqrt{-EG}(\nu_1-\nu_2)={\rm const}\neq 0.\leqno(3.2)$$
\end{prop}

{\it Proof}\,: The equality $\sqrt{-EG}\,(\nu_1-\nu_2)=c\,\lambda \,
\mu,\; c= {\rm const} \neq 0$, and Lemma 3.1 imply the
assertion. \qed
\vskip 2mm

\subsection{Strongly regular time-like W-surfaces.}
We consider strongly regular time-like W-surfaces, i.e. time-like W-surfaces, satisfying the
condition
$$\nu_u(u,v)\nu_v(u,v)\neq 0, \qquad (u,v) \in \mathcal D.$$

Our main theorem for such surfaces is

\begin{thm}\label {T:7.3}
Given two differentiable functions $f(\nu), \, g(\nu); \; \nu \in
\mathcal{I},$ $f(\nu)-g(\nu)\neq 0,$  $f'(\nu)\,g'(\nu)\neq 0$ and a
differentiable function $\nu(u,v), \; (u,v) \in {\mathcal D}$
satisfying the conditions
$$\nu_u\,\nu_v\neq 0, \quad \nu(u,v)\in \mathcal{I}.$$

Let $(u_0, v_0) \in \mathcal D, \; \nu_0=\nu(u_0, v_0)$ and $\mathfrak{a}\neq
0, \,\mathfrak{ b}\neq 0$ be two constants. If
$$\begin{array}{l}
\ds{\mathfrak{a}^2\exp\left(2\int_{\nu_0}^{\nu}
\frac{f'd\nu}{f-g}\right) \left[g'\nu_{uu} +
\left(g''-\frac{2g'^2}{g-f}\right)\nu^2_u\right]}\\
[4mm]
\ds{\;\;\;+\mathfrak{b}^2\exp\left(2\int_{\nu_0}^{\nu}
\frac{g'd\nu}{g-f}\right)\left[f'\nu_{vv}+
\left(f''-\frac{2f'^2}{f-g}\right)\nu^2_v\right]=fg(f-g)},
\end{array}\leqno(3.3)$$
then there exists a unique (up to a motion) strongly regular time-like Weingarten surface\\
$\mathcal M:\; z=z(u,v),$ \, $(u,v)\in \mathcal D_0 \subset \mathcal D$ with
invariants
$$\begin{array}{c}
\nu_1=f(\nu), \quad \nu_2=g(\nu), \\
[2mm]
\displaystyle{\gamma_1=\exp\left(\int_{\nu_0}^{\nu} \frac{g'd\nu}{g-f}\right)\,
\frac{-\mathfrak{b}f'}{f-g}\,\nu_v,
\quad \gamma_2=\exp\left(\int_{\nu_0}^{\nu}
\frac{f'd\nu}{f-g}\right)\, \frac{-\mathfrak{a}g'}{f-g}\,\nu_u.}
\end{array} \leqno(3.4) $$
Furthermore, $(u,v)$ are natural principal parameters for $\mathcal{M}$.
\end{thm}

\emph{Proof}: Taking into account
Proposition \ref {P:5.4}, we obtain that the integrability
conditions 2.1) and 2.2) in Theorem 2.2 reduce to (3.3),
which proves the assertion. \hfill{\qed}
\vskip 2mm

Introducing the functions
$$I:=\int_{\nu_0}^{\nu} \frac{f'(\nu)\,d \nu}{f(\nu)-g(\nu)}\,,\qquad
J:=\int_{\nu_0}^{\nu} \frac{g'(\nu)\,d \nu}{g(\nu)-f(\nu)}\leqno (3.5)$$
we can write the PDE (3.3) in the form
$$\mathfrak{a}^2\,e^{2I}\,\left(J_{uu}+I_u\,J_u-J_u^2 \right)
-\mathfrak{b}^2\,e^{2J}\,\left(I_{vv}+I_v\,J_v-I_v^2\right)=-f\,g,\leqno (3.6)$$
and the principal geodetic curvatures (3.4) in the form
$$\gamma_1=-\mathfrak{b}\,e^J\,I_v,\quad \gamma_2=\mathfrak{a}\,e^I\,J_u.\leqno (3.7)$$

Hence, with respect to natural principal parameters each strongly regular time-like
Weingarten surface possesses a {\it natural PDE} (3.3) (or equivalently (3.6)).

\vskip 2mm
\subsection{Time-like W-surfaces with $\gamma_1=0$.}

In this subsection we consider time-like W-surfaces in Minkowski space with $\gamma_1=0$
and prove the fundamental theorem of Bonnet type for this class.

Let $\mathcal{M}:\; z=z(u,v),\; (u,v)\in \mathcal{D}$ be a  time-like W-surface,
parameterized by natural principal parameters. Then we can assume
$$\mathfrak{a}\,\sqrt E=e^I, \qquad \mathfrak{b}\,\sqrt G=e^J,$$
where $I$ and $J$ are the functions (3.5) and $\mathfrak{a},\,\mathfrak{b}$
are some positive constants.
We note that under the condition $\gamma_1=0$ it follows that the function
$\nu=\nu(u)$ does not depend on $v$.

Considering the system (2.3), we obtain that the compatibility conditions for
this system reduce to only one - the Gauss equation, which has the form:
$$X(\gamma_2)-\gamma_2^2=-f(\nu)\,g(\nu).$$

Thus we obtain the following Bonnet type theorem for time-like W-surfaces
satisfying the condition $\gamma_1=0$:
\begin{thm}
Given two differentiable functions $f(\nu), \; g(\nu); \; \nu \in \mathcal{I},$
$f(\nu)-g(\nu)\neq 0,$ $f'(\nu)\,g'(\nu)\neq 0$ and a differentiable function
$\nu(u,v)=\nu(u), \; (u,v) \in {\mathcal D}$ satisfying the conditions
$$\nu_u\neq 0, \quad \nu(u,v)\in \mathcal{I}.$$

Let $(u_0, v_0) \in \mathcal D, \; \nu_0=\nu(u_0, v_0)$ and $\mathfrak{a}>0$
be a constant. If
$$\mathfrak{a}^2\,e^{2 I}\,(J_{uu}+I_u\,J_u-J_u^2)=-f(\nu)\,g(\nu),\leqno(3.8)$$
then there exists a unique (up to a motion) time-like W-surface
$\mathcal{M}:\; z=z(u,v),$ \, $(u,v)\in \mathcal D_0 \subset \mathcal D$
with invariants
$$\begin{array}{ll}
\nu_1=f(\nu), & \nu_2=g(\nu), \\
[2mm]
\gamma_1=0, & \gamma_2=\mathfrak{a}\,e^{I}\,(J)_u.
\end{array}\leqno(3.9) $$
Furthermore, $(u,v)$ are natural principal parameters on $\mathcal{M}$.
\end{thm}

Hence, with respect to natural principal parameters each time-like
Weingarten surface with $\gamma_1=0$ possesses a {\it natural ODE} (3.8).
\vskip 2mm

\section{Parallel time-like surfaces in Minkowski space and their natural PDE's}

Let $\mathcal{M}: \; z=z(u,v),\; (u,v)\in \mathcal{D}$ be a time-like surface,
parameterized by
principal parameters and $l(u,v),\;l^2=1$ be the unit normal vector field of $\mathcal{M}$.
The parallel surfaces of $\mathcal{M}$ are given by
$$ \overline{{\mathcal{M}}}(a): \; \bar z(u,v)= z(u,v) + a\,l(u,v), \quad a={\rm const} \neq 0,
\quad (u,v)\in \mathcal{D}.\leqno(4.1)$$

We call the family $\{ \overline{{\mathcal{M}}}(a), \;
 a={\rm const} \neq 0\}$ the \emph{parallel
family} of $\mathcal{M}$.

Taking into account (4.1), we find
$$\bar z_{u}=(1-a\,\nu_1)\, z_u,\quad \bar z_{v}=(1-a\,\nu_2)\, z_v.\leqno(4.2)$$

Excluding the points, where $(1-a\,\nu_1)(1-a\,\nu_2)=0$, we obtain that the corresponding
unit normal vector fields $\bar l$ to $ \overline{{\mathcal{M}}}(a)$ and $l$ to $\mathcal{M}$
satisfy the equality
$\bar l=\varepsilon\,l,$ where $\varepsilon:= {\rm sign}\,(1-a\,\nu_1)(1-a\,\nu_2).$
In view of (4.2) it follows that $\bar E< 0$ and $\bar G> 0$.
Hence, the parallel surfaces $ \overline{{\mathcal{M}}}(a)$ of a time-like surface $\mathcal{M}$
are also time-like surfaces.

The relations between the principal curvatures $\nu_1(u,v)$,
$\nu_2(u,v)$ of $\mathcal{M}$ and $\bar \nu_1(u,v)$,  $\bar \nu_2(u,v)$ of its
parallel time-like surface $\overline{{\mathcal{M}}}(a)$ are
$$\bar \nu_1=\varepsilon\,\frac{\nu_1}{1-a\,\nu_1}\,, \quad \bar \nu_2=\varepsilon\,
\frac{\nu_2}{1-a\,\nu_2};\quad \nu_1=\frac{\varepsilon\,\bar\nu_1}{1+a\,\varepsilon\,\bar\nu_1}\,,
\quad \nu_2=\frac{\varepsilon\,\bar\nu_2}{1+a\,\varepsilon\,\bar\nu_2}.\leqno (4.3)$$

Let $K=\nu_1\,\nu_2,\; \displaystyle{H=\frac{1}{2}(\nu_2+\nu_2),\;H'=\frac{1}{2}(\nu_2-\nu_2)}$
be the three invariants of the time-like
surface $\mathcal{M}$. The equalities (4.3) imply the  relations
between the invariants $\bar K$, $\bar H$ and $\bar H'$ of $\overline{{\mathcal{M}}}(a)$
and the corresponding invariants of $\mathcal{M}$:
$$K=\frac{\bar K}{1+2 a\,\varepsilon \bar H + a^2 \bar K}\,, \quad
H=\frac{\varepsilon\,\bar H+a \bar K}{1+2 a\,\varepsilon \bar H + a^2 \bar K}\,,\quad
H'=\frac{\varepsilon\,\bar H'}{1+2 a\,\varepsilon\, \bar H + a^2 \bar
K}\,.\leqno (4.4)$$
\vskip 2mm

Now let $\mathcal{M}: \; z=z(u,v),\; (u,v)\in \mathcal{D}$ be a time-like Weingarten surface
with Weingarten functions $f(\nu)$ and $g(\nu)$. We suppose that $(u,v)$ are
natural principal parameters for $\mathcal{M}$.
We show that $(u,v)$ are also natural principal parameters
for any parallel time-like surface $\overline{{\mathcal{M}}}(a)$.
\begin{prop}
The natural principal parameters $(u,v)$ of a given time-like W-surface $\mathcal{M}$
are natural principal parameters for all parallel time-like surfaces
$ \overline{{\mathcal{M}}}(a),\; a={\rm const} \neq 0$ of $\mathcal{M}$.
\end{prop}

{\it Proof} : Let $(u,v)\in \mathcal{D}$ be natural principal parameters for $\mathcal{M}$,
$(u_0,v_0)$ be a fixed point in $\mathcal D$ and $\nu_0=\nu(u_0,v_0)$.
The coefficients $E$ and $G$ of the first fundamental form of $\mathcal{M}$ are given by (3.1).
The corresponding coefficients $\bar E$ and $\bar G$ of $\overline{ {\mathcal{M}}}(a)$
in view of (4.2) are
$$\bar E=(1-a\,\nu_1)^2\,E,\quad \bar G=(1-a\,\nu_2)^2\,G.\leqno(4.5)$$
Equalities (4.3) imply that $\overline{ {\mathcal{M}}}(a)$ is again a Weingarten
surface with Weingarten functions
$$\bar\nu_1(u,v)=\bar f(\nu)=\frac{\varepsilon f(\nu)}{1-af(\nu)}\,,
\quad\bar\nu_2(u,v)=\bar g(\nu)=\frac {\varepsilon g(\nu)}{1-ag(\nu)}. \leqno(4.6)$$
Using (4.6), we compute
$$\bar f-\bar g=\frac{\varepsilon(f-g)}{(1-a\,f)(1-a\,g)}\,,$$
which shows that ${\rm sign} \,(\bar f-\bar g)={\rm sign}\,(f-g)$.

Further, we denote by $f_0:=f(\nu_0), \; g_0:=g(\nu_0)$ and taking into account
(3.2) and (4.5), we compute
$$\sqrt{-\bar E\,\bar G}\,(\bar f-\bar g)=\sqrt{-E\, G}\,( f- g)=
{\rm const} \neq 0,$$
which proves the assertion. \qed
\vskip 2mm
Using the above statement, we prove the following theorem.

\begin{thm}
The natural PDE of a given time-like W-surface $\mathcal{M}$ is the natural PDE of any parallel
time-like surface $\overline{{\mathcal{M}}}(a),\; a={\rm const}\neq 0$, of $\mathcal{M}$.
\end{thm}
{\it Proof} : We have to express the equation (3.3)  in terms of the
Weingarten functions of the parallel time-like surface $ \overline{{\mathcal{M}}}(a)$.

Putting
$$\bar E_0=(1-a\,\nu_{1}(u_0,v_0))^2\,E_0=-\mathfrak{a}^{-2}\,(1-a\,f_0)^2
=:-\bar {\mathfrak{a}}^{-2},$$
$$\bar G_0=(1-a\,\nu_{2}(u_0,v_0))^2\,G_0=\mathfrak{b}^{-2}\,(1-a\,g_0)^2
=:\bar {\mathfrak{b}}^{-2},$$
we obtain
$$\begin{array}{l}
\ds{\;\;\;\bar{\mathfrak{a}}^2\exp\left(2\int_{\nu_0}^{\nu}
\frac{\bar f'd\nu}{\bar f-\bar g}\right) \left[\bar g'\nu_{uu} +
\left(\bar g''-\frac{2\bar g'^2}{\bar g-\bar f}\right)\nu^2_u\right]}\\
[4mm]
\ds{+\,\bar{\mathfrak{b}}^2\exp\left(2\int_{\nu_0}^{\nu}
\frac{\bar g'd\nu}{\bar g-\bar f}\right)\left[\bar f'\nu_{vv}+
\left(\bar f''-\frac{2\bar f'^2}{\bar f-\bar g}\right)\nu^2_v\right]-
\bar f\,\bar g(\bar f-\bar g)}\\
[4mm]
=\ds{\mathfrak{a}^2\exp\left(2\int_{\nu_0}^{\nu}
\frac{f'd\nu}{f-g}\right) \left[g'\nu_{uu} +
\left(g''-\frac{2g'^2}{g-f}\right)\nu^2_u\right]}\\
[4mm]
\ds{+\,\mathfrak{b}^2\exp\left(2\int_{\nu_0}^{\nu}
\frac{g'd\nu}{g-f}\right)\left[f'\nu_{vv}+
\left(f''-\frac{2f'^2}{f-g}\right)\nu^2_v\right]- f\,g(f-g)}.
\end{array}$$

Hence, the natural PDE of $\overline{{\mathcal{M}}}(a)$ in terms of the Weingarten functions
$\bar f(\nu)$, $\bar g(\nu)$ coinsides with the natural PDE of $\mathcal{M}$ in terms of the
Weingarten functions $f(\nu)$ and $g(\nu)$. \qed
\vskip 2mm

\section{Time-like surfaces whose curvatures satisfy a linear relation}

We now consider time-like W-surfaces, whose three
invariants $K$, $H$ and $H'$ satisfy a linear relation:
$$\delta K=\alpha \,H+\beta\,H'+\gamma, \quad \alpha, \beta, \gamma, \delta - {\rm constants},
\quad \alpha^2-\beta^2+4\gamma\delta\neq 0.\leqno(5.1)$$
\vskip 1mm

A time-like W-surface with principal curvatures $\nu_1$ and
$\nu_2$ is said to be {\it linear fractional} if
$$\nu_1=\frac{A\nu_2+B}{C\nu_2+D}\,, \quad BC-AD\neq 0.\leqno(5.2)$$

We exclude the case $A=D,\,B=C=0$, which characterizes the points with $H^2-K=0$
and show that the classes of surfaces with characterizing conditions (5.1) and (5.2),
respectively, coincide.

\begin{lem}
Any surface whose invariants
$\displaystyle{K=\nu_1\,\nu_2,\; H=\frac{1}{2}(\nu_1+\nu_2),\; H'=\frac{1}{2}(\nu_1-\nu_2)}$
satisfy a linear relation $(5.1)$
is a linear fractional time-like Weingarten surface determined by $(5.2)$, and vice versa.
\end{lem}

The relations between the constants $\alpha, \beta, \gamma, \delta$ in (5.1) and $A, B, C, D$
in (5.2) are given by the equalities:
$$\alpha=A-D, \quad \beta = -(A+D), \quad \gamma=B,
\quad \delta=C.\leqno(5.3)$$

We denote by $\mathfrak{K}$ the class of all time-like surfaces with $H^2-K>0$, whose
curvatures satisfy (5.1) or equivalently (5.2).

The aim of our study is to classify all natural PDE's of the surfaces from the class
$\mathfrak{K}$.

The parallelism between two surfaces
given by (4.1) is an equivalence relation. On the other hand, Theorem 4.2 shows that
the surfaces from an equivalence class have one and the same natural PDE. Hence, it is
sufficient to find the natural PDE's of the equivalence classes. For any equivalence class,
we use a special representative, which we call \emph{a basic class}. Thus the classification
of the natural PDE's of the surfaces in the class $\mathfrak{K}$ reduces to the
natural PDE's of the basic classes.

In view of Theorem 4.2, we prove the following classification theorem.
\begin{thm}
Up to similarity, the time-like surfaces in Minkowski space, whose
curvatures $K$, $H$ and $H'$ satisfy the linear relation
$$\delta K = \alpha H + \beta H' + \gamma, \quad \alpha, \beta, \gamma, \delta -
{\rm constants}; \quad \alpha^2-\beta^2 + 4 \gamma\delta \neq 0,$$
are described by the natural PDE's of the following basic surfaces:
\vskip 3mm

$(1)$
$H=0:\quad  \nu=e^{\lambda},\quad \bar\Delta \lambda =  e^{\lambda};$\\
\vskip 3mm

$(2)$
$H=\displaystyle{\frac{1}{2}:\quad \nu=\frac{1}{2}(1-e^{\lambda}),\quad
\bar\Delta \lambda = \sinh{\lambda}};$\\
\vskip 3mm

$(3)$
$H'=1: \quad \bar\Delta^*(e^{\nu})=2\,\nu\,(\nu+2);$\\
\vskip 3mm

$(4)$
$H=\beta\,H'\;(\beta^2>1):\qquad \bar\Delta^*(\nu^{\beta})=
\displaystyle{2\,\frac{\beta\,(\beta+1)}{(\beta-1)^2}\,\nu;}$\\
\vskip 3mm

$(5)$
$H=\beta\,H'\; (\beta^2<1):\qquad \Delta^*(\nu^{\beta})=
\displaystyle{2\,\frac{\beta\,(\beta+1)}{(\beta-1)^2}\,\nu;}$\\
\vskip 3mm

$(6)$
$\left|\begin{array}{l}
H=\beta\,H'+1 \\
[1mm]
\beta^2>1\end{array}\right.:\quad
 \displaystyle{\nu=\frac{(\beta-1)\,\lambda+2}{2},\quad
\bar\Delta^*(\lambda^{\beta})= \frac{\beta\,((\beta-1)\lambda+2)((\beta+1)\lambda+2)}
{2\,(\beta-1)\,\lambda};}$\\
\vskip 3mm

$(7)$
$\left|\begin{array}{l}
H=\beta\,H'+1 \\
[1mm]
\beta^2<1\end{array}\right.:\quad
 \displaystyle{\nu=\frac{(\beta-1)\,\lambda+2}{2},\quad
\Delta^*(\lambda^{\beta})= \frac{\beta\,((\beta-1)\lambda+2)((\beta+1)\lambda+2)}
{2\,(\beta-1)\,\lambda};}$\\
\vskip 3mm

$(8)$
$K=-1:\quad \nu=\tan\lambda,\quad \Delta \lambda = -\sin \lambda;$\\
\vskip 3mm

$(9)$
$K=2\,H': \quad \displaystyle{\nu=\frac{\lambda-4}{\lambda-2},\quad
\bar\Delta^*(e^{\lambda})=2;}$\\
\vskip 3mm

$(10)$
$K=\beta\,H'+\gamma \; (\beta\neq 0,\,\gamma<0):\quad
 \left | \begin{array}{l}
\displaystyle{\nu=\lambda+\frac{\beta}{2}},\quad
\displaystyle{\mathcal{I}=\frac{1}{\sqrt{-\gamma}}\,\arctan\frac{\lambda}{\sqrt{-\gamma}}}\,,\\
[4mm]
 \displaystyle{\bar\Delta^*(e^{\beta\,\mathcal{I}}) =-\frac{\beta\,\gamma}{2}\,
\frac{\lambda\,\left(\beta\,\lambda+2\,\gamma\right)}{\lambda^2-\gamma}.}
\end{array} \right.$
\end{thm}
\vskip 2mm

{\it Proof} :
According to the constant $C$ in (5.2), the linear fractional time-like W-surfaces
are divided into two classes:
linear fractional time-like W-surfaces, determined by the condition $C=0$
and linear fractional time-like W-surfaces, determined by
the condition $C\neq 0$.
\vskip 2mm

{\bf I. Linear fractional time-like Weingarten surfaces with $C=0$.}

This class is determined by the equality
$$\alpha\,H+\beta\,H'+\gamma=0,\quad (\alpha,\gamma)\neq (0,0),\quad
\alpha^2-\beta^2\neq 0.\leqno(5.4)$$

For the invariants of the time-like parallel surface $\overline{{\mathcal M}}(a)$ of
$\mathcal M$, because of (4.4),
we get the relation
$$\varepsilon\,(\alpha+2\,a\,\gamma)\,\bar H+\varepsilon\,\beta\,\bar H'+\gamma
=-a\,(\alpha+a\,\gamma)\,\bar K.\leqno (5.5)$$

Let $\;\eta:= {\rm sign}\,(\alpha^2-\beta^2).$
Each time choosing appropriate values for the constants $\mathfrak{a}$, $\mathfrak{b}$
and $\nu_0$ in (3.3), we consider the following subclasses and their
natural PDE's:
\vskip 2mm

\begin{itemize}
\item[1)] $\alpha=0,\;\beta\neq 0,\;\gamma\neq 0$. Assuming that $\gamma=1$, the relation
(5.4) becomes
$$\beta\,H'+1=0. $$

The natural PDE for these W-surfaces is
$$(e^{-\beta\,\nu})_{uu}-(e^{\beta\,\nu})_{vv}=
\frac{2}{\beta}\,\nu\,(\beta\,\nu-2).\leqno (5.6)$$

Up to similarities these time-like W-surfaces are generated by the basic class $H'=1$
with the natural PDE
$$(e^{\nu})_{uu}-(e^{-\nu})_{vv}=2\,\nu\,(\nu+2),\leqno (5.6^*)$$
which is the case (3) in the statement of the theorem.
\vskip 2mm

\item[2)]
 $\displaystyle{\alpha\neq 0,\;\gamma=0}$. Assuming that $\alpha=1$, the relation
 (5.4) becomes
$$H+\beta\,H'=0.$$
\begin{itemize}
\item[2.1)] $\beta\neq 0,\;\eta=-1\; (\beta^2-1>0).$
Choosing
$\displaystyle{\mathfrak{b}^2\,\frac{\beta-1}{\beta+1}\,\nu_0^{-(\beta+1)}=1,\;
\mathfrak{a}^2\,\nu_0^{\beta-1}=1}$, the natural PDE becomes
$$\left(\nu^{-\beta}\right)_{uu}-\left(\nu^{\beta}\right)_{vv}=
2\,\frac{\beta(\beta-1)}{(\beta+1)^2}\,\nu,\leqno(5.7)$$
which is the case (4) in the statement of the theorem.
\vskip 2mm

\item[2.2)] $\beta \neq 0,\;\eta=1\; (\beta^2-1<0)$. Choosing
$\displaystyle{\mathfrak{b}^2\,\frac{\beta-1}{\beta+1}\,\nu_0^{-(\beta+1)}=-1,\;
\mathfrak{a}^2\,\nu_0^{\beta-1}=1}$, the natural PDE becomes
$$\left(\nu^{-\beta}\right)_{uu}+\left(\nu^{\beta}\right)_{vv}=
2\,\frac{\beta(\beta-1)}{(\beta+1)^2}\,\nu,\leqno (5.8)$$
which is the case (5) in the statement of the theorem.
\vskip 2mm

\item[2.3)] $\beta = 0$. Putting  $\;\nu =e^{\lambda}$, we get the natural
PDE for time-like surfaces with $H=0$:
$$\lambda_{uu}-\lambda_{vv} = e^{\lambda},\leqno (5.9)$$
which is the case (1) in the statement of the theorem.
\end{itemize}
\vskip 2mm

\item[3)] $\alpha\neq 0,\;\beta=0,\; \gamma\neq 0.$ Assuming that $\alpha=1$,
the relation (5.4) becomes
$$H+\gamma=0.$$
Putting  $\;\displaystyle{| H|\,e^{\lambda}:=H-\nu}=H'>0$,
we get the one-parameter system of natural PDE's for CMC time-like surfaces
with $H=-\gamma$:
$$\lambda_{uu}-\lambda_{vv}=2\,|H|\,\sinh \lambda.\leqno(5.10)$$

Up to similarities these time-like W-surfaces are generated by the basic class $|H|=\frac{1}{2}$
with the natural PDE
$$ \lambda_{uu}-\lambda_{vv}= \sinh \lambda,\leqno(5.10^*)$$
which is the case (2) in the statement of the theorem.
\vskip 2mm

\item[4)] $ \alpha\neq 0,\;\beta\neq 0,\;\gamma\neq 0.$

Assuming that $\alpha=1$ we have
$$H+\beta\,H'+\gamma=0,\quad \beta^2-1\neq 0.$$
Let $\displaystyle{\lambda:=2\,H'=\frac{-2}{\beta+1}\,(\nu+\gamma)}>0$.

\begin{itemize}
\item[4.1)] If $\eta=-1 \; (\beta^2-1>0)$ and choosing
$$\mathfrak{b}^2=\frac{\beta+1}{\beta-1}\,\left(\frac{-2}{\beta+1}
(\nu_0+\gamma)\right)^{\beta+1},\;
\mathfrak{a}^2\,=\left(\frac{-2}{\beta+1}(\nu_0+\gamma)\right)^{-(\beta-1)},$$

the natural PDE becomes
$$\left(\lambda^{-\beta}\right)_{uu}-\left(\lambda^{\beta}\right)_{vv}=
\frac{\beta}{2\,(\beta+1)}\,
\frac{((\beta+1)\lambda+2\,\gamma)((\beta-1)\lambda+2\,\gamma)}{\lambda}\,.
\leqno (5.11)$$

Up to similarities these time-like W-surfaces are generated by the basic class\\
$H=\beta\,H'+1,\;\beta^2>1$ with the natural PDE
$$\left(\lambda^{\beta}\right)_{uu}
-\left(\lambda^{-\beta}\right)_{vv}=
\frac{\beta}{2\,(\beta-1)}\,
\frac{((\beta+1)\lambda+2)((\beta-1)\lambda+2)}{\lambda}\,,
\leqno (5.11^*)$$
which is the case (6) in the statement of the theorem.
\vskip 2mm

\item[4.2)] If $\eta=1\; (\beta^2-1<0)$ and choosing
$$\mathfrak{b}^2=-\frac{\beta+1}{\beta-1}\,\left(\frac{-2}{\beta+1}
(\nu_0 + \gamma)\right)^{\beta+1},\;
\mathfrak{a}^2\,=\left(\frac{-2}{\beta+1}(\nu_0 + \gamma)\right)^{-(\beta-1)},$$

the natural PDE becomes
$$\left(\lambda^{-\beta}\right)_{uu}+\left(\lambda^{\beta}\right)_{vv}=
\frac{\beta}{2\,(\beta+1)}\,
\frac{((\beta+1)\lambda+2\,\gamma)((\beta-1)\lambda+2\,\gamma)}{\lambda}\,.
\leqno (5.12)$$

Up to similarities these time-like W-surfaces are generated by the basic class\\
$H=\beta\,H'+1,\;\beta^2<1$ with the natural PDE
$$\left(\lambda^{\beta}\right)_{uu}
+\left(\lambda^{-\beta}\right)_{vv}=
\frac{\beta}{2\,(\beta-1)}\,
\frac{((\beta+1)\lambda+2)((\beta-1)\lambda+2)}{\lambda}\,,
\leqno (5.12^*)$$
which is the case (7) in the statement of the theorem.
\end{itemize}
\end{itemize}

{\bf II. Linear fractional time-like Weingarten surfaces with $C\neq 0$.}

Let $C=1$. The equality (5.1) gets the form
$$K=\alpha\,H+\beta\,H'+\gamma.\leqno(5.13)$$

The corresponding relation for the parallel surface $ \overline{{\mathcal M}}(a)$ is
$$\varepsilon(\alpha+2\,a\,\gamma)\,\bar H+\varepsilon\,\beta\,\bar H'+\gamma=
(1-a\,\alpha-a^2\,\gamma)\,\bar K.\leqno(5.14)$$

Each time choosing appropriate values for the constants
$\mathfrak{a}$, $\mathfrak{b}$ and $\nu_0$ in (3.3), we consider the following subclasses
and their natural PDE's:
\begin{itemize}
\item[5)] $\alpha = \gamma= 0, \; \beta \neq 0$. The relation (5.13) becomes
$$K=\beta H' \quad \iff \quad \rho_1-\rho_2=-\frac{2}{\beta}\,,$$
where $\displaystyle{\rho_1=\frac{1}{\nu_1},\;\rho_2=\frac{1}{\nu_2}}$ are the principal
radii of curvature of $ \mathcal{M}$.

Putting $\lambda:=\displaystyle{4\,\frac{\nu-\beta}{2\,\nu-\beta}}$,
the natural PDE of these time-like surfaces gets the form
$$\left(e^\lambda\right)_{uu} -\left(e^{-\lambda}\right)_{vv}-\frac{\beta^4}{8}=0.
\leqno (5.15)$$

Up to similarities these time-like W-surfaces are generated by the basic class
$K=2\,H'$ with the natural PDE
$$\left(e^\lambda\right)_{uu} -\left(e^{-\lambda}\right)_{vv}-2=0,
\leqno (5.15^*)$$
which is the case (9) in the statement of the theorem.
\vskip 2mm

\item [6)] $(\alpha, \gamma)\neq (0,0), \; \alpha^2+4\gamma\geq 0$. The relation (5.14)
implies that there exists a time-like surface $\overline{{\mathcal M}}(a)$,
parallel to $\mathcal M$, which satisfies the
relation (5.4). Hence the natural PDE of $ \mathcal{M}$ is one of the PDE's (5.6) - (5.12).
\vskip 2mm

\item[7)] $\alpha^2+4\,\gamma <0$. It follows that $\gamma <0$. The relation (5.14)
implies that there exists a time-like surface $\overline{{\mathcal M}}(a)$ parallel
to $\mathcal M$, which satisfies the relation
$$K=\beta H'+\gamma.\leqno (5.16)$$
\begin{itemize}
\item[7.1)] $\beta = 0$. The relation (5.16) becomes
$K=\gamma<0,$ i.e. $\overline{{\mathcal M}}$ is of constant negative sectional curvature
$\gamma$.
Putting $\displaystyle{\lambda:=2\,\arctan\frac{\nu}{\sqrt{-\gamma}}}$,
we get the natural PDE of this surface
$$ \lambda_{uu}+\lambda_{vv}=-K^2\,\sin \lambda.
\leqno (5.17)$$

Up to similarities these time-like W-surfaces are generated by the basic class
$K=-1$ with the natural PDE
$$\lambda_{uu}+\lambda_{vv}=-\sin \lambda,\leqno (5.17^*)$$
which is the case (8) in the statement of the theorem.
\vskip 2mm

\item[7.2)] $\beta \neq 0$, $\gamma< 0$. Choosing
$\displaystyle{\nu_0=\frac{\beta}{2}}$,
 the natural PDE of $\overline{\mathcal{M}}$ becomes
 $$(\exp{(\beta\,\mathcal{I})})_{uu}-
(\exp{(-\beta\,\mathcal{I})})_{vv}
=-\frac{\beta\,\gamma}{2}\,
\frac{\lambda\,\left(\beta\,\lambda+2\,\gamma\right)}{\lambda^2-\gamma}\,,\leqno(5.18)$$
where
$$\mathcal{I}=\frac{1}{\sqrt{-\gamma}}
\,\arctan\frac{\lambda}{\sqrt{-\gamma}},\quad
\lambda:=\nu-\frac{\beta}{2},$$
which is the case (10) in the statement of the theorem.
\end{itemize}
\end{itemize}
\hfill{$\square$}
\vskip 4mm

\section{Summary}

Summarizing the results in \cite{GM2, GM3} and in the present paper,
we obtain the following parallel between the natural PDE's
describing linear fractional W-surfaces in $\R^3$, linear fractional
space-like and time-like W-surfaces in $\R^3_1$, respectively.
\begin{itemize}
\item[(i)] The natural PDE for a Weingarten surface in Euclidean space is of the type:
$$\begin{array}{l}
\ds{\mathfrak{a}^2\exp\left(2\int_{\nu_0}^{\nu}
\frac{f'd\nu}{f-g}\right) \left[g'\nu_{uu} +
\left(g''-\frac{2g'^2}{g-f}\right)\nu^2_u\right]}\\
[4mm]
-\ds{\;\;\;\mathfrak{b}^2\exp\left(2\int_{\nu_0}^{\nu}
\frac{g'd\nu}{g-f}\right)\left[f'\nu_{vv}+
\left(f''-\frac{2f'^2}{f-g}\right)\nu^2_v\right]}=-fg(f-g),
\end{array}$$
or equivalently
$$\mathfrak{a}^2\,e^{2 I}\,(J_{uu}+I_u\,J_u-J_u^2)+
\mathfrak{b}^2\,e^{2 J}(I_{vv}+I_v\,J_v-I_v^2)
=f(\nu)\,g(\nu).$$
\vskip 4mm

\item[(ii)] The natural PDE for a space-like Weingarten surface in Minkowski space is
of the type:
$$\begin{array}{l}
\ds{a^2\exp\left(2\int_{\nu_0}^{\nu}
\frac{f'd\nu}{f-g}\right) \left[g'\nu_{uu} +
\left(g''-\frac{2g'^2}{g-f}\right)\nu^2_u\right]}\\
[4mm]
\ds{\;\;\;-b^2\exp\left(2\int_{\nu_0}^{\nu}
\frac{g'd\nu}{g-f}\right)\left[f'\nu_{vv}+
\left(f''-\frac{2f'^2}{f-g}\right)\nu^2_v\right]=fg(f-g)},
\end{array}$$
or equivalently
$$\mathfrak{a}^2\,e^{2 I}\,(J_{uu}+I_u\,J_u-J_u^2)+
\mathfrak{b}^2\,e^{2 J}(I_{vv}+I_v\,J_v-I_v^2)
=-f(\nu)\,g(\nu).$$
\vskip 4mm

\item[(iii)] The natural PDE for a time-like Weingarten surface with real principal
curvatures in Minkowski space is of the type:
$$\begin{array}{l}
\ds{\mathfrak{a}^2\exp\left(2\int_{\nu_0}^{\nu}
\frac{f'd\nu}{f-g}\right) \left[g'\nu_{uu} +
\left(g''-\frac{2g'^2}{g-f}\right)\nu^2_u\right]}\\
[4mm]
\ds{\;\;\;+\mathfrak{b}^2\exp\left(2\int_{\nu_0}^{\nu}
\frac{g'd\nu}{g-f}\right)\left[f'\nu_{vv}+
\left(f''-\frac{2f'^2}{f-g}\right)\nu^2_v\right]=fg(f-g)},
\end{array}$$
or equivalently
$$\mathfrak{a}^2\,e^{2I}\,\left(J_{uu}+I_u\,J_u-J_u^2 \right)
-\mathfrak{b}^2\,e^{2J}\,\left(I_{vv}+I_v\,J_v-I_v^2\right)=-f(\nu)\,g(\nu).$$
\end{itemize}
\vskip 2mm

Therefore for the corresponding basic linear fractional surfaces in $\mathbb{R}^3$
and $\mathbb{R}^3_1$ we obtain the correspondence between their natural PDE's.

\end{document}